\documentclass[12pt, a4paper]{article}
\usepackage[left=1in,right=1in,top=1.0in,bottom=1in]{geometry}
\usepackage{amsmath,amssymb,amsthm,amsfonts,afterpage,hyperref}
\usepackage{amscd,subeqnarray}
\usepackage{graphicx}
\usepackage{caption}
\usepackage{subcaption}
\usepackage{authblk}
\usepackage{tikz}
\usepackage{array}
\usepackage{wrapfig}
\usepackage{color}
\usepackage{ragged2e}
\usepackage{stmaryrd}
\usepackage{siunitx}
\usepackage{multirow}
\usepackage{xcolor,cancel}
\usepackage{boldfonts}
\usepackage{symbols}
\usepackage{footmisc}
\usepackage{float}
\usepackage{setspace}
\usepackage{bm}
\usepackage{booktabs}
\usepackage{tikz, xcolor}
\usepackage{soul}
\usepackage[T1]{fontenc}
\usepackage[utf8]{inputenc}
\usepackage[linesnumbered,ruled,vlined]{algorithm2e}
\SetKwInput{KwInput}{Input}                
\SetKwInput{KwOutput}{Output} 
\usepackage{cite}
\usetikzlibrary{shapes,arrows, positioning,calc}	
\usetikzlibrary{arrows,snakes,backgrounds}
\newcommand{\bfa}[1]{\boldsymbol{#1}} 			
			
\newcommand{\bfeps}{\boldsymbol{\epsilon}}

\newcommand{\Sym}{\text{Sym}}   			%
\newcommand{\curl}{\text{curl}}   				%
\newcommand{\tr}{\text{tr}}       				%

\DeclareMathAlphabet{\mathpzc}{OT1}{pzc}{m}{it}

\newcommand{\bfu}{\boldsymbol{u}}	
	
\newcommand{\bfv}{\boldsymbol{v}}	

\newcommand{\bfF}{\boldsymbol{F}}	
\newcommand{\bfE}{\boldsymbol{E}}	
\newcommand{\bfC}{\boldsymbol{C}}
\newcommand{\bfB}{\boldsymbol{B}}	
\newcommand{\bfx}{\boldsymbol{x}}	
\newcommand{\bfX}{\boldsymbol{X}}	
	
\newcommand{\bfT}{\boldsymbol{T}}		
\newcommand{\bfI}{\boldsymbol{I}}	 

\newtheorem{theorem}{Theorem}

\newtheorem{remark}{Remark}

\newtheorem*{cf}{Strong formulation}
\newtheorem*{cwf}{Continuous weak formulation}

\newtheorem*{dfep}{Discrete Finite Element Problem}

\providecommand{\keywords}[1]
{
  \small	
  \textbf{\textit{Keywords---}} #1
}

\title{An efficient finite element method for computing the response of a strain-limiting elastic solid containing a v-notch and inclusions}

\author[1]{Shylaja G.\thanks{g\_shylaja@blr.amrita.edu}}
\author[1]{Kesavulu Naidu V.\thanks{v\_kesavulu@blr.amrita.edu}\thanks{corresponding author}}
\author[1]{Venkatesh B.  \thanks{b\_venkatesh@blr.amrita.edu}\thanks{corresponding author}}
\author[2]{S. M. Mallikarjunaiah\thanks{m.muddamallappa@tamucc.edu}\thanks{corresponding author}}
\affil[1]{Department of Mathematics, Amrita School of Engineering, Amrita Vishwa Vidyapeetham, Bengaluru, 560035, INDIA}
\affil[2]{Department of Mathematics \& Statistics,
Texas A\&M University-Corpus Christi,
Corpus Christi, TX 78412, USA}
\date{}

\begin{document}
\maketitle
\begin{abstract}
Accurate triangulation of the domain plays a pivotal role in computing the numerical approximation of the differential operators. A good triangulation is the one which aids in reducing discretization errors.  In a standard collocation technique, the smooth curved domain is typically triangulated with a mesh by taking points on the boundary to approximate them by polygons. However, such an approach often leads to geometrical errors which directly affect the accuracy of the numerical approximation. To restrict such geometrical errors, \textit{isoparametric}, \textit{subparametric}, and \textit{iso-geometric} methods were introduced which allow the approximation of the curved surfaces (or curved line segments). In this paper, we present an efficient finite element method to approximate the solution to the elliptic boundary value problem (BVP), which governs the response of an elastic solid containing a v-notch and inclusions. The algebraically nonlinear constitutive equation along with the balance of linear momentum reduces to second-order quasi-linear elliptic partial differential equation.  Our approach allows us to represent the complex curved boundaries by smooth \textit{one-of-its-kind} point transformation. The main idea is to obtain higher-order shape functions which enable us to accurately compute the entries in the finite element matrices and vectors. A Picard-type linearization is utilized to handle the nonlinearities in the governing differential equation. The numerical results for the test cases show considerable improvement in the accuracy.  

\end{abstract}
\noindent \keywords{ Finite element method; Curved elements; Point transformation; Strain-limiting elastic body; Quasi-linear eliptic partial differential equation.} 

\section{Introduction}

The real-world problems are often described by linear or nonlinear partial differential equations posed on complex geometric domains \cite{naidu2013advantages,sasikala2023efficient}.  Simulating such scenarios possesses a lot of challenges arising from nonlinearity in the differential operators \cite{vasilyeva2023,gou2023computational,gou2023finite}, geometric complications \cite{MVSMM2023}, nonlocality in the boundary conditions \cite{ferguson2015numerical}. In many instances, the differential operators are approximated by any of the standard collocation methods such as \textit{finite difference methods} or \textit{finite element methods} or \textit{boundary element methods} or \textit{spectral methods} or even in some cases \textit{meshless methods}. Creating a good computational mesh for a better approximation in the aforementioned methods is very challenging. In clinical research, the images rendered by computer-aided design software must be accurately represented by meshes (or triangulation). For the finite element computations, the triangulation with straight sides (for example, as in triangles or quadrilaterals) offers a straightforward way of computing the two-dimensional or three-dimensional numerical integrations using high-accurate quadrature rules. This is mainly due to the availability of the standard shape functions (like Lagrange type or any other standard ones) as well as access to the quadrature rules on the regular polygons. Such is not the case when the domain of interest has curved boundaries. The triangulation of such a domain needs to have triangles (or quadrilaterals) with one side curved and the other straight sides to accurately represent the computational domain so that the discretization is minimal (or even zero). The issue is very important when the curved inclusions or heterogeneities are present in the region interesting for computations. Hence, a school of thought of creating curved elements was introduced in \cite{ergatoudis1968curved,mcleod1975use,mitchell1979advantages,zienkiewicz2005finite}. For such curved finite element methods (CFEM), the approximation of the differential operators directly depends on the veracity of the geometrical representation \cite{ciarlet1972interpolation,scott1973finite}. In CFEM, there are three main approaches: first, the entire global domain containing curved boundary parts is transformed into a standard polygonal shape, and a global computing method can be achieved; second, triangulate the curved region into the regular domain; third, isoparametric finite elements. There are some inherent challenges in these methods \cite{rathod2008use,nagaraja2010use,naidu2010use,naidu2013advantages}. In this paper, we explore the parabolic arcs in matching the curved boundaries (which are the circular inclusions in the domain) by \textit{unique} point transformations and obtain a higher-order finite element method for cubic order curved triangular element. We explore the numerical technique for approximating the solution of a quasi-linear partial differential equation. The overall algorithm involves Picard's linearization to handle the nonlinearities and a special finite element method for spatial discretization.

 The mathematical model studied in this work governs the response of a limiting-strain elastic solid formulated within the new theory developed in \cite{rajagopal2007elasticity,rajagopal2003implicit,rajagopal2011non,rajagopal2014nonlinear,rajagopal2009class,rajagopal2011modeling,gou2015modeling,MalliPhD2015,bulivcek2015analysis,bulivcek2015existence,bulivcek2014elastic}. Such a setup offers an attractive framework to model the response of large class of materials under both mechanical and thermal loading \cite{yoon2022finite,bonito2020finite}. It is well documented that some materials behave nonlinearly even in less than $2\%$ strain, hence the necessity of developing first-order, yet nonlinear constitutive relations.  More importantly, a numerical tool is much needed to correctly approximate the curved regions in the computational domain. Hence, this paper is the first attempt in such a direction. 
 
 This paper is organized as follows: In section~\ref{intr_imp_theory} we introduce the implicit constitutive theory for the response of the elastic solid. A second quasi-linear elliptic partial differential equation is obtained for the anit-plane shear loading setup is also presented in section~\ref{intr_imp_theory}. A BVP and a theorem for the existence of a solution are given in section~\ref{bvp_existence}. Section~\ref{fem_disc} contains a finite element discretization of the BVP. The construction of the shape functions, discrete formulation along with the computational algorithm are all given in the same section. The numerical examples with domains containing v-notch and inclusions are provided in section~\ref{num_exp}. The conclusions and some information about the future works are given in section~\ref{conclusion}. 

\section{Implicit theory of elasticity} \label{intr_imp_theory}
In this section, our goal is to provide a brief introduction to an implicit nonlinear theory formulated in a series of papers by Rajagopal 
\cite{rajagopal2003implicit,rajagopal2007elasticity,rajagopal2011non,rajagopal2011conspectus,rajagopal2014nonlinear,rajagopal2007response} for the response of the bulk material.  The constitutive relationships are algebraically nonlinear, however, using such relations we develop a mathematical model for the response of an isotropic linear elastic body containing both v-notch and inclusions. The ultimate aim of this investigation is to provide a convergent numerical method for the discretization of the BVP. 

Let $\mathcal{B}$ denote an elastic body that is in equilibrium with the externally applied mechanical loading and we suppose that $\partial\mathcal{B}$ is the boundary. The elastic material is assumed to occupy a region in two-dimensional space.  Let $\bfX$  and $\bfx$ are the points in the reference and deformed configurations of the body.  Let ``$\Sym$'' be the linear space of symmetric tensors equipped with standard \textit{Frobenius norm}.  We suppose $\bfu \colon \mathcal{B} \to \mathbb{R}^2$ denote the displacement, $\bfT \colon \mathcal{B} \to \mathbb{R}^{2 \times 2}_{\Sym}$ denote the Cauchy stress, and $\bfeps \colon \mathcal{B} \to \mathbb{R}^{2 \times 2}_{\Sym}$ denotes the infinitesimal strain tensor. The linear space of symmetric tensors equipped with standard \textit{Frobenius norm}.  Let $\bfF \colon \mathcal{B} \to \mathbb{R}^{2 \times 2}$ denote the deformation gradient,  $\bfC \colon \mathcal{B} \to \mathbb{R}^{2 \times 2}$ denote the right Cauchy-Green stretch tensor, $\bfB \colon \mathcal{B} \to \mathbb{R}^{2 \times 2}$ denote the left Cauchy-Green stretch tensor, $\bfE \colon \mathcal{B} \to \mathbb{R}^{2 \times 2}$ denote the Lagrange strain, respectively defined by:
\begin{subequations}
\begin{align}
\bfF &:=  \bfI + \nabla \bfu, \;
 \bfB :=\bfF\bfF^{\mathrm{T}}, \; \bfC :=\bfF^{\mathrm{T}}\bfF,  \; \bfE := \dfrac{1}{2} \left( \bfC - \bfI \right),\\
\bfeps &:= \dfrac{1}{2} \left( \nabla \bfu + \nabla \bfu^{\mathrm{T}}\right),
\end{align}
\end{subequations}
where $\left( \cdot \right)^{\mathrm{T}}$ denotes the \textit{transpose} operator, $\bfI$ is the identity tensor, respectively. The infinitesimal strain theory assumes that 
\begin{equation}\label{small_grad}
\max_{\bfX \in \mathcal{B}} \| \nabla_X \bfu \|  \ll \mathcal{O}(\delta), \quad \delta \ll 1.
\end{equation}
The above assumption \eqref{small_grad} implies
\begin{equation}\label{lin_results}
\bfB \approx \bfI + 2 \bfeps, \;\;  \bfC \approx \bfI + 2 \bfeps , \;\;  \bfE \approx  \bfeps, \;\; \det \bfF = 1 + \tr(\bfeps)
\end{equation}
Rajagopal \cite{rajagopal2007elasticity} generalized the constitutive relationships of Cauchy elasticity by introducing  implicit elastic response relations as 
\begin{equation}\label{implicit-1}
0 = \mathcal{F}(\bfT,  \;  \bfB ).
\end{equation}
where $\mathcal{F}$ is tensor-valued, the \textit{isotropic function} such that it obeys \cite{rajagopal2007elasticity,rajagopal2011modeling,gou2015modeling,mallikarjunaiah2015direct,MalliPhD2015}
\begin{equation}
\mathcal{F}(\boldsymbol{Q}\bfT\boldsymbol{Q}^{\mathrm{T}}, \; \boldsymbol{Q}\bfB\boldsymbol{Q}^{\mathrm{T}}) = \boldsymbol{Q}  \mathcal{F}(\bfT,  \;  \bfB )\boldsymbol{Q}^{\mathrm{T}} \quad \forall \boldsymbol{Q} \in \mathbb{O},
\end{equation} 
where $\mathbb{O}$ indicates the orthogonal group. Rajagopal \cite{rajagopal2007elasticity} considered a special subclass of \eqref{implicit-1}:
\begin{equation}\label{SL1}
\bfB := \mathcal{F}( \bfT), \quad \mbox{with}  \quad \sup_{\bfT \in \Sym}   \| \mathcal{F}( \bfT) \|  \leq M, \;\; M >0. 
\end{equation}
If such a constant $M$ exists,  then the response relations \eqref{SL1} are referred to as \textit{strain-limiting} \cite{MalliPhD2015,mallikarjunaiah2015direct}. Applying the standard linearization procedure \eqref{small_grad} to \eqref{SL1},  we obtain 
\begin{subequations}
\begin{align}
\bfeps  &=\mathcal{F}( \bfT), \\
&= \beta_1 \, \bfI + \beta_2 \, \bfT + \beta_2 \, \bfT^2,
\end{align}
\end{subequations}
Finally, the governing system of equations to model the behavior of an elastic material within the framework of algebraically nonlinear   theories is given by
\begin{subequations}
\begin{align}
-\nabla \cdot \bfT &=\bf{0}, \quad \mbox{and} \quad \bfT = \bfT^T, \label{equilib:eq} \\
\bfeps &= \Psi_{0}\left( \tr \,  \bfT, \; \| \bfT \|  \right) \bfI + \Psi_{1}\left(  \| \bfT \|  \right) \bfT, \quad {\Psi_{0}\left( 0, \; \cdot  \right)=0}, \label{eqn:main} \\
\curl \, \curl \, \bfeps &=\bfa{0}, \label{eqn_scompata} \\ 
\bfeps &:= \frac{1}{2} \left( \nabla \bfu +  \nabla \bfu^{T}  \right). \label{eqn_linstrain}
\end{align}
\end{subequations}
Where both $\Psi_{0} \colon \mathbb{R} \times \mathbb{R}_{+} \to \mathbb{R}$ and $\Psi_{1} \colon \mathbb{R}_{+} \to \mathbb{R}$ are functions of invariants of the Cauchy stress. 

\begin{remark}
In the above governing system of equations, equation~\eqref{equilib:eq}$_1$ is from the balance of linear momentum for the quasi-static situation and with the body force acting on the material and equation~\eqref{equilib:eq}$_2$ is from the balance of angular momentum. Equation~\eqref{eqn:main} is the constitutive relationship governing the response of the material to mechanical/thermal stimuli. Equation~\eqref{eqn_scompata} is the strain-compatibility condition and equation~\eqref{eqn_linstrain} is the equation for the linearized strain. 
\end{remark}

\subsection{Anti-plane Shear Problem (Mode-III)}
Our main objective in this work is to provide an efficient method to simulate the response of a strain-limiting elastic body containing v-notch and inclusions and to the authors' best knowledge, this is the first attempt to investigate a simulation tool for such a setup of algebraically nonlinear constitutive relations. 
The problem is a static v-notch in a material that is under mode-III loading. In the current loading situation the displacement $\bfa{u}$ is a scalar function of $x$ and $y$, i.e. 
\begin{equation}\label{eq:disp_vector}
\bfa{u}(x,\, y) = \left( 0, 0, u( x,\, y)  \right).
\end{equation}
 The constitutive relation for the linear, isotropic, and homogeneous material is 
 \begin{equation}
\bfT = 2 \, \mu \, \bfeps,
\end{equation}
where $\mu$ is the shear modulus (N/m$^2$) of the material. Then the constitutive equation within \textit{strain-limiting theory of elasticity} is 
\begin{equation}
\bfeps = \Psi_{1}\left( \| \bfT \|  \right) \bfT. \label{eqn:main2}
\end{equation}
Let us introduce \textit{Airy's stress function} $\Phi$ and the stress components through the gradient of $\Phi$ as
\begin{equation}
\bfT_{13} = \dfrac{\partial \Phi}{\partial y}, \quad \bfT_{23} = - \, \dfrac{\partial \Phi}{\partial x},
\end{equation}
which readily satisfies the balance of linear momentum 
\begin{equation}\label{blm}
\dfrac{\partial \bfT_{13}}{\partial x} + \dfrac{\partial \bfT_{23}}{\partial y}=0. 
\end{equation}  
The strain-compatibility equation \eqref{eqn_scompata} reduces to the system of equations for the strain components:
\begin{subequations}
\begin{align}
\dfrac{\partial }{\partial x} \left( \dfrac{\partial }{\partial y}  \bfeps_{13}  - \dfrac{\partial }{\partial x}  \bfeps_{23}   \right) &=0, \\
\dfrac{\partial }{\partial y} \left( \dfrac{\partial }{\partial x}  \bfeps_{23}  - \dfrac{\partial }{\partial y}  \bfeps_{13}   \right) &=0. 
\end{align}
\end{subequations}
Notice that the components of the strain tensor depend only on the in-plane coordinates $x$ and $y$, we obtain that 
\begin{equation}\label{strain_compatibility}
 \dfrac{\partial }{\partial y}  \bfeps_{13}  - \dfrac{\partial }{\partial x}  \bfeps_{23}  =0. 
\end{equation}
From equation~\eqref{eqn:main}, the components of the strain tensor be written as:
\begin{equation}\label{strain_components}
\bfeps_{13}= \Psi_{1}\left(  \| \bfT \|  \right) \bfT_{13}, \quad \bfeps_{23}= \Psi_{1}\left(  \| \bfT \|  \right) \bfT_{23}. 
\end{equation}
Notice that we have used the simpler constitutive relationship with only the second term in \eqref{eqn:main}. Then using \eqref{strain_components} in \eqref{strain_compatibility} and also utilizing the Airy's stress function $\Phi$ via the balance of linear momentum \eqref{blm}, we obtain a second-order quasilinear elliptic partial differential equation for $\Phi$ 
\begin{equation}\label{pde:nlin1}
  - \nabla \cdot \left( \Psi_{1} \left( \| \nabla \Phi \| \right) \;  \nabla \Phi \right) =0,
\end{equation}
with 
\begin{equation}\label{eqn:grad_airy_norm}
\|\nabla \Phi \| = \sqrt{ \left( \dfrac{\partial \Phi}{\partial x} \right)^2 +  \left( \dfrac{\partial \Phi}{\partial y} \right)^2}.
\end{equation}
For the development of a BVP, we use the following form for $\Psi_{1}$:  
\begin{equation}\label{eq:Psi}
 \Psi_{1} (\|\bfT\|) = \frac{1}{2 \, \mu \left( 1 + \, \| \bfT \| \right)}.
\end{equation}
\begin{remark}
It is clear that $\Psi_1$ is both invertible and monotone and generates hyperelastic relations and nonhyperelastic in some cases \cite{mai2015strong,mai2015monotonicity,rajagopal2007response,rajagopal2011conspectus,rajagopal2004thermomechanical}. 
It is clear that the monotonically decreasing function $\Psi_{1} (r)$ satisfies 
\begin{equation}\label{EqLimit}
\lim_{r \to \infty} r \Psi_{1} (r) \to  \widehat{c},
\end{equation}
where $\widehat{c}$ depends on the material parameters. 
An important consequence of \eqref{EqLimit} is when $\| \bfT \| \to \infty$ the ``strains'' will be uniformly bounded in the entire material body including at the vicinity of concentrators such as crack-tip or re-entrant corner. 
\end{remark}
 In the view of Equation~\eqref{eq:Psi}, the nonlinear PDE (Equation~\eqref{pde:nlin1}) now takes the form
\begin{equation} \label{pde:mech}
- \nabla \cdot \left( \frac{\nabla \Phi}{2 \, \mu \left( 1 + \; \|\nabla \Phi \| \; \right)}  \right) = 0.
\end{equation}
\begin{remark}
The second-order quasi-linear elliptic partial differential equation derived in this paper is similar to the one studied in  \cite{bulivcek2015existence,bulivcek2014elastic,bulivcek2015analysis} and also similar to the minimal surface equation from the calculus of variations. Proof for the existence and uniqueness of weak solutions is provided in \cite{bulivcek2015existence} for a v-notch boundary value problem. 
\end{remark}
\begin{remark}
There are several experimental studies available in the literature  \cite{hao2005super,saito2003multifunctional} which show the clear nonlinear stress-strain behavior for several materials 
such as gum metal \cite{kulvait2019state}, titanium alloys \cite{tian2015nonlinear,devendiran2017thermodynamically}, alloys such as Ti-30Nb-10Ta-5Zr (TNTZ-30). The nonlinear behavior is observed well within $1\%$ strain.  Such behavior is traditionally modeled using a linear relationship between the ``infinitesimal strain'' and Cauchy stress using the classical linearized elasticity construct, which is odd as a linear relationship can't fit the nonlinear data. The aforementioned nonlinear behavior can be modeled within the implicit theory of elasticity \cite{rajagopal2014nonlinear,muliana2018determining,kowalczyk2019finite}.  Rocks \cite{bustamante2020novel} and many rubber-like materials \cite{bustamante2021new} can be modeled using a special subclass of models within the larger class of strain-limiting nonlinear relationships. In \cite{bustamante2021new},  a constitutive model for rubber is proposed using the principle stresses as the main variables, wherein a better corroboration with experimental data is shown than Ogden's model \cite{ogden1972large}. The new material models studied by several aforementioned researchers are good candidates to fit the experimental data on gum metals, rubber-like materials, and high-strength titanium alloys  \cite{rajagopal2014nonlinear,devendiran2017thermodynamically,kulvait2019state,bustamante2021new}. 
\end{remark}

\begin{figure}[H]
     \centering
     \begin{subfigure}[b]{0.4\textwidth}
         \centering
         \includegraphics[width=\textwidth]{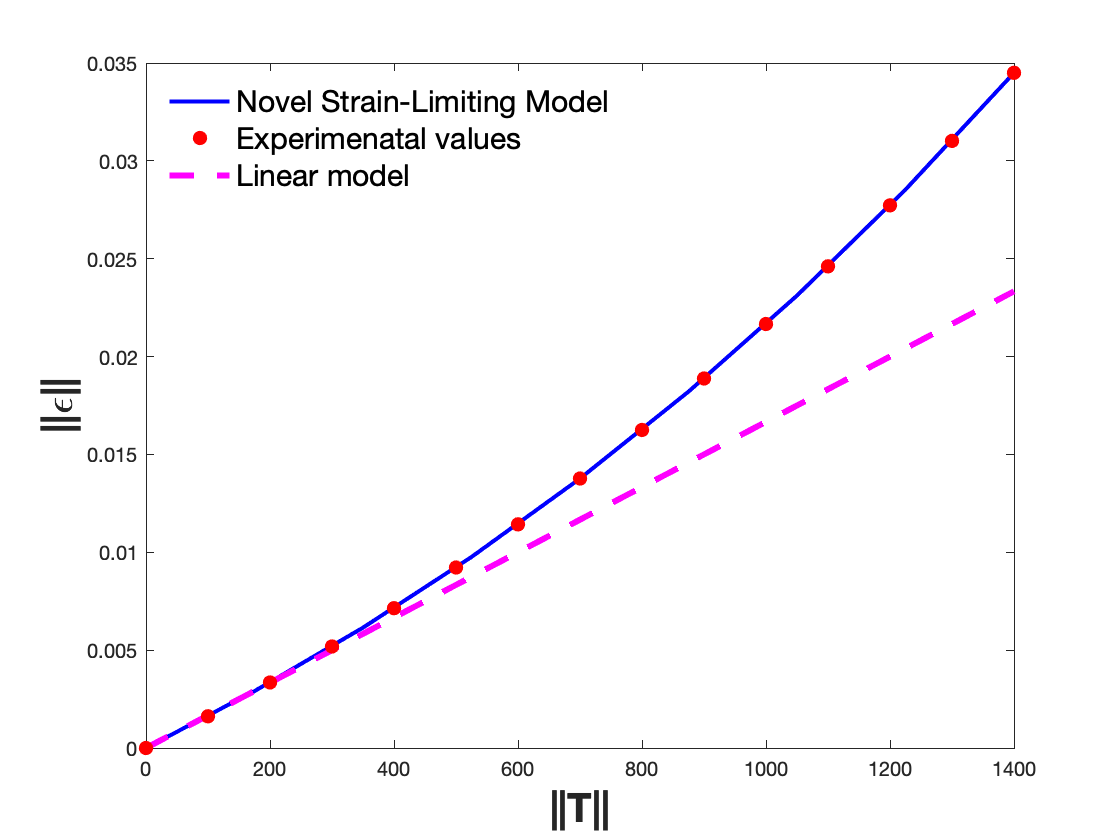}
         \caption{Gum metal experimental data \cite{kulvait2017modeling}}
     \end{subfigure}
     \hspace{1cm}
     \begin{subfigure}[b]{0.4\textwidth}
         \centering
         \includegraphics[width=\textwidth]{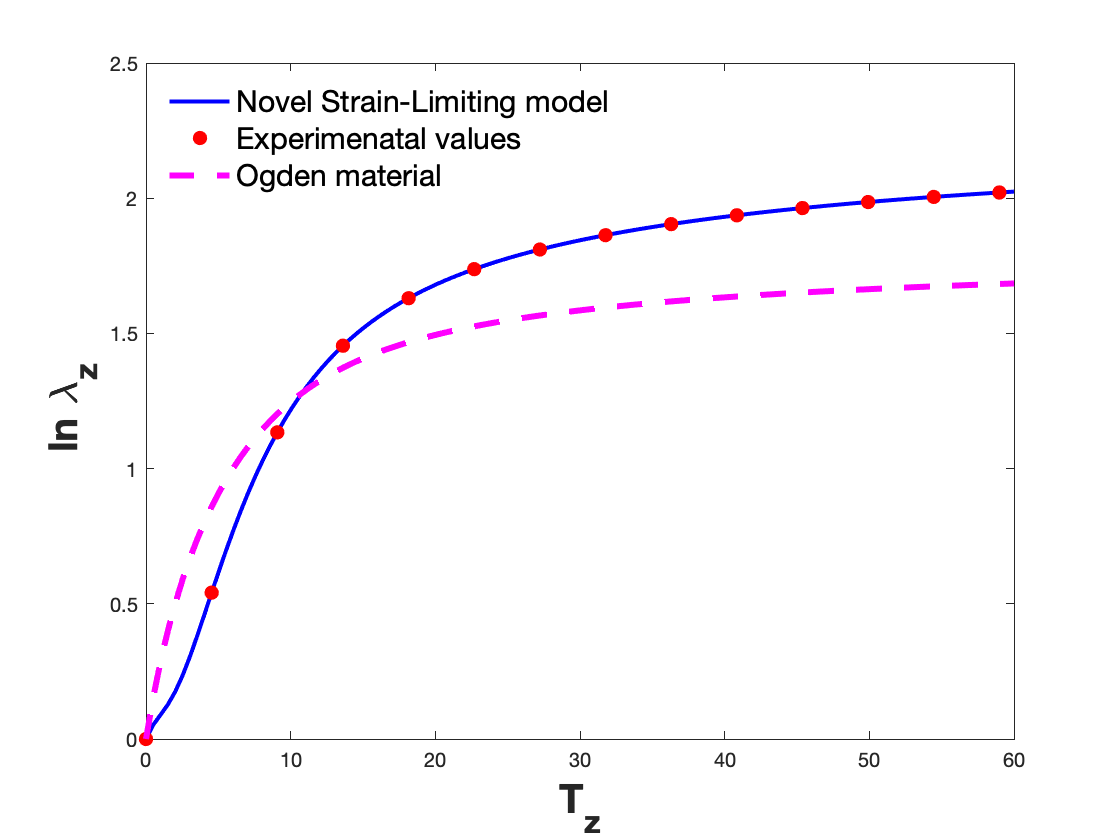}
         \caption{ Behavior of rubber-like solids \cite{ogden1972large}}
     \end{subfigure}
     \caption{Plot of the best fit for the model similar to the one in \eqref{eqn:main2}, for a specific choice of material invariants, compared with the experimental data for gum metal and rubber-like material.}
              \label{exp_plot}
\end{figure}
It is clear from the above Figure-\ref{exp_plot} that for some experiments, the classical linear model and the \textit{Ogden material model} \cite{ogden1972large} can't best fit the nonlinear behavior in the ``\textit{small strain}'' regime (less than $1\%$ strains). Therefore, there is a great opportunity for modelers and computational scientists to develop new constitutive relationships that can correctly describe the experimentally observed phenomenon in many high-strength materials used in the construction industry, orthopedic implants, and aerospace applications.  

\section{Boundary value problem and existence of solution}\label{bvp_existence}
We assume that the response of the body is constituted through the algebraically nonlinear relationship whereas the material is still geometrically linear. We suppose that the solid body is homogeneous and initially unstressed. Let $\mathcal{B} = \Omega \times \mathbb{R}$ where $\Omega \subset  \mathbb{R}^2$ is a simply connected domain with its boundary $\partial \Omega$ being $C^{0, \, 1}$. The boundary consists of two non-overlapping parts $\Gamma_D$ and $\Gamma_N$. Both of these boundary parts may contain finitely many parts and one can parameterize each one of them. Let $L^{p}(\Omega)$ be the space of all \textit{Lebesgue integrable functions} with $p\in[1, \infty)$ and $\left( \cdot, \; \cdot \right)_{L^p}$ denotes the inner product and $\| \cdot \|_{L^2}$ be the corresponding induced norm. The  classical \textit{Sobolev space}  is denoted by $H^{1}(\Omega)$ with the norm defined through:
\[
 \| v \|_{H^{1}(\Omega)} := \| v \|^2 + \| \nabla v \|^2, 
 \]
and let $H^{-1}(\Omega)$ denote the dual space to $H^{1}(\Omega)$. We also define the following subspaces of $H^1(\Omega)$
\begin{equation}
V := H^1(\Omega), \quad V^{0} := \left\{ v \in  H^{1}(\Omega) \colon v=0 \; \mbox{on} \; \Gamma_D\right\}. 
\end{equation}
The out-of-plane displacement function $w(x, \; y)$ is the only unknown variable that defines the mechanics of the material body.  For the loading setup on hand, we have only two non-zero components in both stress and strain tensors. The problem under consideration is simplified to
\begin{cf}
Given the material parameter $\mu$, find $\Phi \in C^2(\Omega)$ such that 
\begin{subequations}\label{f1a}
\begin{align}
-  \, \nabla \cdot \left( \frac{\nabla \Phi}{2 \, \mu \left( 1 + \; \|\nabla \Phi \|\; \right)}  \right)  &=0, \; \mbox{in} \;\; \Omega, \;\;  \label{f1_1} \\
 \Phi &= \widehat{\Phi}, \quad  \mbox{on} \quad \partial \Omega  \label{f1_2}
\end{align}
\end{subequations}
\end{cf}
In the above formulation, we have turned an original traction and displacement problem into a pure non-homogeneous Dirichlet-type problem for one single function $\Phi$. Our formulation is similar to the one studied in \cite{bulivcek2014elastic,bulivcek2015existence,kulvait2019state,kulvait2013anti,HyunMSM_MMS2022}. 

A well-posed weak formulation is derived from the above strong formulation by multiplying \eqref{f1_1} by test function $\varphi \in V^0$ and integrating by parts over the domain $\Omega$. 
\begin{cwf}
Given the material parameter $\mu$, find $\Phi \in V$ such that 
\begin{equation}\label{wf}
 \left( \frac{\nabla \Phi}{2 \, \mu \left( 1 + \; \|\nabla \Phi \| \; \right)}, \;\; \nabla \varphi  \right)  =0, \quad \forall \varphi \in V^0 \\
\end{equation}
\end{cwf}

The above continuous weak formulation has been shown to possess a unique solution. Here we recast the main theorem in \cite{bulivcek2015existence} to gain a useful result about the existence and uniqueness of the solution to our formulation. We have the following theorem from  \cite{bulivcek2015existence}:

\begin{theorem}
$\Omega \subset \mathbb{R}^2$ is a simply connected domain with $C^{0,1}$ boundary consisting of two subsets  $\partial \Omega_1 \cup \partial \Omega_2$ so that $\partial \Omega_1 := \left\{ \bfa{x} \in \partial \Omega \colon \bfa{t} (\bfa{x}) \cap \Omega = \emptyset        \right\}$, where $\bfa{t} (\bfa{x})$ is the tangent to the boundary. Both the boundaries may consist of finitely many parts and each can be parameterized to an open interval. Each part in  $\partial \Omega_2$ is flat. There exists a function $\Phi^0 \in W^{1,\, \infty}$ such that $\Phi^0  \big|_{\partial \Omega} = \widehat{\Phi}$. Let $\alpha \in (0, \, 2)$, then there is a unique solution $\Phi \in W^{1, \, 2}_{loc}(\Omega)$ for the problem satisfying
\begin{subequations}
\begin{align}
\Phi - \Phi^0  &\in W^{1, \, 1}_{0}(\Omega) \\
\left( \frac{\nabla \Phi}{2 \, \mu \left( 1 + \; \|\nabla \Phi \|^{\alpha} \; \right)^{1/\alpha}}, \;\; \nabla \psi  \right)  &=0, \quad \forall \psi \in W^{1, \, 1}_{0}(\Omega)
\end{align}
\end{subequations}  
\end{theorem}

The above theorem is the same as the one proved in  \cite{bulivcek2015existence}  has two additional modeling parameters $\beta$ and $\alpha$, however, for $\beta=1$ and $\alpha=1$ our continuous formulation falls well within the purview of the proof achieved in   \cite{bulivcek2015existence}. Hence our continuous weak formulation in \eqref{wf} is well-posed with the existence of a unique solution.  \\

The above strong formulation doesn't possess any exact analytical solution, even in the simpe setting of 1D. Therefore we look at developing a stable convergent numerical method using finite element discretization to approximate the solution to the weak formulation  \eqref{wf}. In the next section, we propose a stable finite-element discretization in which the shape functions are calculated by using a one-of-its-kind point transformation. The exact transformation and the integration procedure produce superior numerical results. 

\section{Finite element discretization}\label{fem_disc}
In this contribution, we propose a finite element discretization of a nonlinear BVP that models the response of a strain-limiting elastic material containing both v-notch and inclusions. For the numerical approximation, the challenges are the mesh discretization and the PDE approximation. In the case of triangulation of the domain, we utilize the curved triangular elements which do not contain discretization error since the exact geometry is utilized in creating the finite element mesh.  For the numerical approximation,  cubic-order shape functions are adapted and these are obtained from one-of-its-kind point transformation established in the literature \cite{mcleod1975use,rathod2008use}. The effective use of the higher order curved triangular elements are clearly demonstrated in works of \cite{shylaja2021improved, SHYLAJA2021twodim, Murali2019DBF1, Shylaja2019poisson, Murali2019DB, Murali2019DBF2}  An exact analytical solution is not known for the physical setup considered in this paper.  Hence we compute a relative difference of the numerical solution between two successive Picard iterations as a measure to quantify the quality of the approximation proposed. 

Let $\Omega$ be a two-dimensional open domain with $\Gamma$ being its boundary which is assumed to be sufficiently smooth.  The boundary $\Gamma$  is composed of the sets $\Gamma_{\Phi}$ and $\Gamma_h^i$ such $\Gamma = \Gamma_{\Phi} \cup_{i=1}^{m} \Gamma_h^i$ with $m \in \mathbb{N}$ being the number of inclusions. Here $\Gamma_{\Phi}$ is the part where we apply the Dirichlet boundary condition and $\Gamma_h^i$ are the internal inclusions. The internal inclusions $\Gamma_h^i$ are assumed to be traction-free. The region near the holes is triangulated using ``curved'' triangles which contain one-side curved and two-side straight lines. 

For approximating the  displacement field  $\bfu$,  we define the following space,
\begin{equation}
S_h = \left\{  \Phi_h \in  \left( C(\overline{\Omega})\right) \colon \left. \Phi_h\right|_K \in \mathbb{Q}_3, \; \forall K \in \mathcal{T}_h \right\},
\end{equation}
where $\mathbb{Q}_3$ is a set of tensor-products of polynomials up to an order of $2$ over the reference cell $\widehat{K}$. Then the discrete  approximation space is
\begin{equation}\label{app-spaces}
\widehat{V}_h = S_h \, \cap  \, V. 
\end{equation}

\subsection{Construction of shape functions}
\begin{figure}[h!]
 \centering
 \noindent\begin{subfigure}[b]{0.4\textwidth}
  \includegraphics[width=1\textwidth]{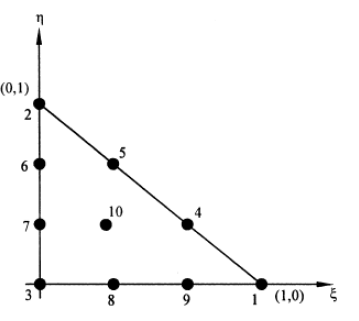}
  \caption{Unmapped triangular element of cubic order}
 \end{subfigure}%
 \qquad
\noindent\begin{subfigure}[b]{0.4\textwidth}
  \includegraphics[width=1\textwidth]{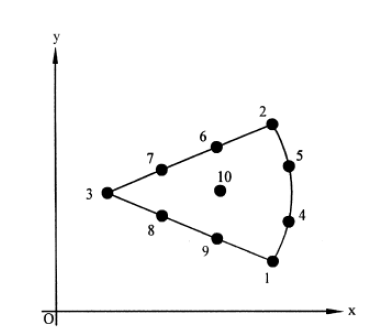}
  \caption{Mapped triangular element of cubic order}
 \end{subfigure}%
 \caption{Mapping the 10-node cubic order with one side curved triangle to the standard right-angled unit triangle}
 \label{fig_mapping}
\end{figure}

Consider the triangular elements in which one of the sides is curved and the other two sides are straight as shown in Fig. \ref{fig_mapping}. The Lagrange interpolants for the field variable $\Phi$ governing the physical problem are
\begin{equation}\label{eq_phi}
\Phi_h=\sum\limits_{i=1}^\frac{(n+1)(n+2)}{2}{ N_{i}^{(n)}(\xi ,\eta )\Phi_{i}^{e}}  \quad n=3
\end{equation}
where $N_{i}^{(n)}(\xi ,\eta )$ refers to the conventional triangular element shape functions Eq.(\ref{SF}) as derived in \cite{mcleod1975use,rathod2008use} of order $n$ at the node $i$. Hence the transformation formulae between the physical (Cartesian) and the local (natural) coordinate systems are
\begin{equation}\label{eq_subpara}
t=\sum\limits_{i=1}^\frac{(n+1)(n+2)}{2}{ N_{i}^{(n)}(\xi ,\eta )t_{i}}
\end{equation}
The nodes along the straight sides $3-1$ and $3-2$ in Fig.\ref{fig_mapping} are equispaced. Now if we use the standard formulae dividing a line segment in a given ratio from the plane analytical geometry to the straight sides $3-1$ and $3-2$ then the Eq.(\ref{eq_subpara}) reduces to
\begin{equation}
    t(\xi,\eta)=t_3+(t_1-t_3)\xi+(t_2-t_3)\eta +a_{11}^{(n)}(t)\xi\eta+H(n-3)\sum\limits_{{i+j=n, i\neq j}}a_{ij}^{(n)}\xi^i\eta^j, \quad (1\leq i, j \leq n-1),
\end{equation}
where $t$ is the nodal value of the triangular element and $H(n-3)$ is the well-known Heaviside step function or unit step function and it has the meaning for the present as,
\begin{equation}
 H(n-3) = \begin{cases}
                 0 \quad n<3\\
                 1 \quad n \leq 3
 \end{cases}
\end{equation}

The following are the shape functions for a curved triangle element of cubic order:
\begin{subequations}\label{SF}
\begin{align}
N_{1}^{(3)}(\xi ,\,\,\eta ) &=\frac{9{{\xi }^{3}}}{2}-\frac{9{{\xi }^{2}}}{2}+\xi,  \\
N_{2}^{(3)}(\xi ,\,\,\eta ) &=\frac{9{{\eta }^{3}}}{2}-\frac{9{{\eta }^{2}}}{2}+\eta,  \\
N_{3}^{(3)}(\xi ,\,\,\eta ) &=-\frac{9{{\xi }^{3}}}{2}-\frac{9{{\eta }^{3}}}{2}-\frac{27{{\xi }^{2}}\eta }{2}-\frac{27\xi {{\eta }^{2}}}{2}+9{{\xi }^{2}}+9{{\eta }^{2}}+18\xi \eta -\frac{11\xi }{2}-\frac{11\eta }{2}+1,  \\
N_{4}^{(3)}\,(\xi ,\,\eta ) &=\frac{27{{\xi }^{2}}\eta }{2}-\frac{9\xi \eta }{2} \\
N_{5}^{(3)}\,(\xi ,\,\eta ) &=\frac{27\xi {{\eta }^{2}}}{2}-\frac{9\xi \eta }{2} \\
N_{6}^{(3)}\,(\xi ,\,\eta ) &=-\frac{27{{\eta }^{3}}}{2}-\frac{27\xi {{\eta }^{2}}}{2}+18{{\eta }^{2}}+\frac{9\xi \eta }{2}-\frac{9\eta }{2} \\
N_{7}^{(3)}\,(\xi ,\,\eta ) &=\frac{27{{\eta }^{3}}}{2}+27\xi {{\eta }^{2}}-\frac{45{{\eta }^{2}}}{2}+\frac{27{{\xi }^{2}}\eta }{2}-\frac{45\xi \eta }{2}+9\eta, \\
N_{8}^{(3)}(\xi ,\,\,\eta ) &=\frac{27{{\xi }^{3}}}{2}+27{{\xi }^{2}}\eta +\frac{27\xi {{\eta }^{2}}}{2}-\frac{45{{\xi }^{2}}}{2}-\frac{45\xi \eta }{2}+9\xi, \\
N_{9}^{(3)}\,(\xi ,\,\eta ) &=-\frac{27{{\xi }^{2}}}{2}-\frac{27{{\xi }^{2}}\eta }{2}+18{{\xi }^{2}}+\frac{9\xi \eta }{2}-\frac{9\xi }{2}, \\
N_{10}^{(3)\,}(\xi ,\,\eta ) &=-27\xi {{\eta }^{2}}-27{{\xi }^{2}}\eta +27\xi \eta, 
\end{align}
\end{subequations}
Figure \ref{fig_mapping} depicts the domain with two straight sides and one curved side. The subparametric transformations given below are used to convert the global coordinates of a typical triangle element in cubic order to the local coordinates. The cubic order's point transformation is as follows:
\begin{equation}
    t(\xi ,\,\eta )= {{t}_{3}}+({{t}_{1}}-{{t}_{3}})\xi +({{t}_{2}}-{{t}_{3}})\,\eta+\frac{9}{4}[({{t}_{4}}+{{t}_{5}})-({{t}_{1}}+{{t}_{2}})]\xi \eta ,\quad  t=(x,y), \quad ,   n= 3 
\end{equation}
with the boundary node $t_{5}$ as,
\begin{equation}
t_{5}=t_{4}-\frac{1}{3} (t_{1} -t_{2} ), \\
\end{equation}
and the interior node  $t_{10}$ as,
\begin{equation}
t_{10} =\frac{1}{12} (t_{1} +t_{2} +4t_{3} +3t_{4} +3t_{5} )
\end{equation} 

\subsection{Discrete formulation and Picard's iteration} 
Starting with the continuous formulation, the discrete weak formulation is obtained using the shape functions from the approximation space defined in \eqref{app-spaces}. However, the difficulty is in handling the nonlinearities in the continuous case. These can be handled either at the differential equation level or at the linear algebra level. We have used the former method by the well-known \textit{Picard's type linearization}. Such a linearization produces a sequence of linear problems. For fast convergence of the above iterative algorithm, a proper initial guess $\Phi^{0}_h$ is required. In our implementation, we solved the linear problem first (obtained from \eqref{discrete-wf} without the denominator) and subsequently used the computed solution as an initial guess for Picard's iterations.  All our numerical simulations converged within a reasonable number of iterations. Overall, our formulation possesses a unique solution at the discrete level.  Finally, the discrete finite element problem reads: 
\begin{dfep}
Given the Dirichlet boundary data $\Phi^{0}_h \in \widehat{V}_h$, and the $n^{th}$ iteration solution $\Phi^n_h \in \widehat{V}_h$, for $n=0, 1, 2, \cdots $, find $\Phi^{n+1}_h  \in \widehat{V}_h$ such that 
\begin{equation}\label{discrete-wf}
   a(\Phi_h^n; \, \Phi^{n+1}_h,\, v_h) = l(\bfv_h), \forall\, v_h \in \widehat{V}_h,  
\end{equation}
where the bilinear and linear terms are given by
\begin{subequations}\label{A-L-Def}
\begin{align}
a(\Phi_h^n; \, \Phi^{n+1}_h,\, v_h) &=\int_{\Omega} \left[  \frac{\nabla \Phi_h^{n+1}  \cdot  \nabla v_h} {1 +  \;  \| \nabla  \Phi_h^{n} \|}         \right]  \;   d\bfx\, ,  \label{disc_A} \\
l(\bfv_h) &= \int_{\Gamma} {f} \cdot v_h \; d\bfx + \int_{\Gamma} {g} \cdot v_h \; ds.  \label{disc_L}  
\end{align}
\end{subequations}
\end{dfep}
In the above weak formulation, the term $f$ is the known function and $g$ is the applied traction. \\

The following algorithm depicts the overall discrete finite element computational procedure to obtain the numerical solution of the BVP.
 
\begin{algorithm}[H]
\SetAlgoLined
\KwInput{Choose the parameters: $I_{Max}$ (maximum of the iteration number), $Tol$}
Start with a sufficiently refined mesh; \\
Solve the linear problem to obtain the initial guess $\Phi^{0}_h \in \widehat{V}_h$; \\
\For {$n=0, \,1, \, 2, \, \ldots$}{
 \While{$[\text{Iteration Number} < \text{$I_{Max}$}]$ AND $[\text{Residual} > \text{$Tol$}]$}{
  Assemble Equations~\eqref{disc_A} and \eqref{disc_L} using test functions from the discrete finite element space $S_h$ \;
 Use a \textit{direct solver} to solve for $ \Phi_h^{n+1}$\;
  }
  }  
 Save the final converged solution $\Phi_h$ to output files for the visualization\;
 \caption{Algorithm for the finite element computations}
 \label{algo001}
\end{algorithm}


\section{Numerical Examples}\label{num_exp}
To illustrate the efficiency of the proposed method, we consider three examples: the first one consists of nonlinear BVP with a manufactured solution; the second one is the nonlinear BVP on the v-notch domain; the third one is the nonlinear BVP on the v-notch domain containing several inclusions. The challenge in the third example is the meshing near the curved inclusions. The proposed method takes very few elements to discretize the region near inclusions, more importantly, there is no discretization error due to considering curved triangles.

\subsection{Square domain with a manufactured solution}
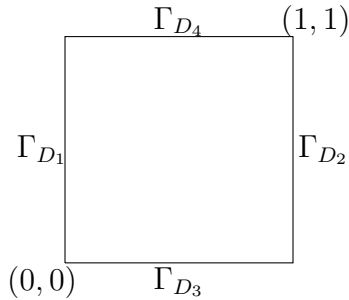
\begin{figure}[H]
\centering
\begin{tikzpicture}
\draw (0,0) -- (3,0) -- (3,3) -- (0,3) -- (0,0);

\node at (-0.3,-0.25)   {$(0,0)$};
\node at (3.3,3.2)   {$(1,1)$};

\node at (-0.3, 1.5)   {$\Gamma_{D_{1}}$};
\node at (3.4, 1.5)   {$\Gamma_{D_{2}}$};
\node at (1.5, -0.25)   {$\Gamma_{D_{3}}$};
\node at (1.5, 3.2)   {$\Gamma_{D_{4}}$};
\end{tikzpicture}
\caption{A domain and the Dirichlet boundary conditions for h-convergence study.}
\label{fig_domain}
\end{figure}

In this example, we consider the nonlinear BVP posed on a square domain $\Omega$ and we choose $\Phi = \frac{\pi}{2} y^2$ as the exact solution. The right-hand side term $f$ is computed by plugging $\Phi$ into the strong form. The boundary values are set by using the manufactured solution $\Phi$. Such a situation allows us to do a convergence analysis in the standard $l_2$ measure. The computational domain is depicted in Figure~\ref{fig_domain}. The boundary conditions are shown in the Table~\ref{ex1_bc}.

\begin{table}[H]
\centering
\caption{Boundary conditions for square domain}
\begin{tabular}{p{0.9 in} p{0.5 in}}\hline
 Boundary & Values \\ \hline
 $\Gamma_{D_{1}}$ & $\frac{\pi}{2} y^2$ \\
 $\Gamma_{D_{2}}$ & $\frac{\pi}{2} y^2$ \\
 $\Gamma_{D_{3}}$ & 0 \\
 $\Gamma_{D_{4}}$ &  $\frac{\pi}{2}$  \\
\end{tabular}
 \label{ex1_bc}
\end{table}

For the illustrative purpose, we have considered three different types of meshes for computations. The mesh information is given in the following table:
\begin{table}[H]
\centering
\caption{Details of computational mesh}
\begin{tabular}{p{0.6 in} p{0.5 in} p{1.25 in}}\hline
Elements & DOF & Boundary Nodes\\ \hline
8 & 49 & 24 \\
16 & 85 & 24\\
32 & 169 & 48\\ \hline
\end{tabular}
\label{tab_square}
\end{table}

\begin{figure}[H]
 \centering
 \noindent\begin{subfigure}[b]{0.25\textwidth}
  \includegraphics[width=1\textwidth]{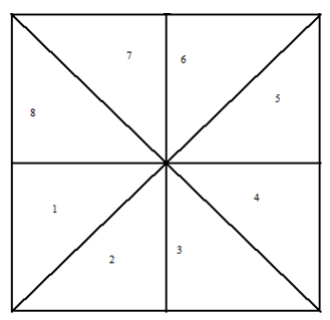}
  \caption{8 elements}
 \end{subfigure}%
 \qquad
\noindent\begin{subfigure}[b]{0.25\textwidth}
  \includegraphics[width=1\textwidth]{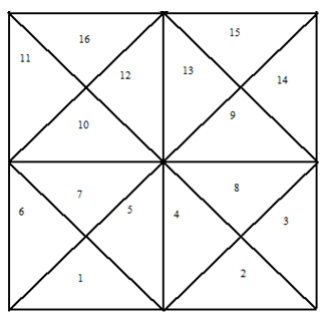}
  \caption{16 elements}
 \end{subfigure}%
 \qquad
\noindent\begin{subfigure}[b]{0.25\textwidth}
  \includegraphics[width=1\textwidth]{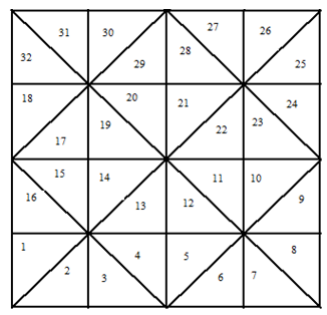}
  \caption{32 elements}
 \end{subfigure}%
 \caption{Doamin discretization into 8, 16, and 32 elements}
 \label{fig_dis_square}
\end{figure}

The discretization of the mesh using different numbers of elements is given in Figure~\ref{fig_dis_square}. 

\begin{figure}[H]
 \centering
\noindent\begin{subfigure}[b]{0.5\textwidth}
  \includegraphics[width=1\textwidth]{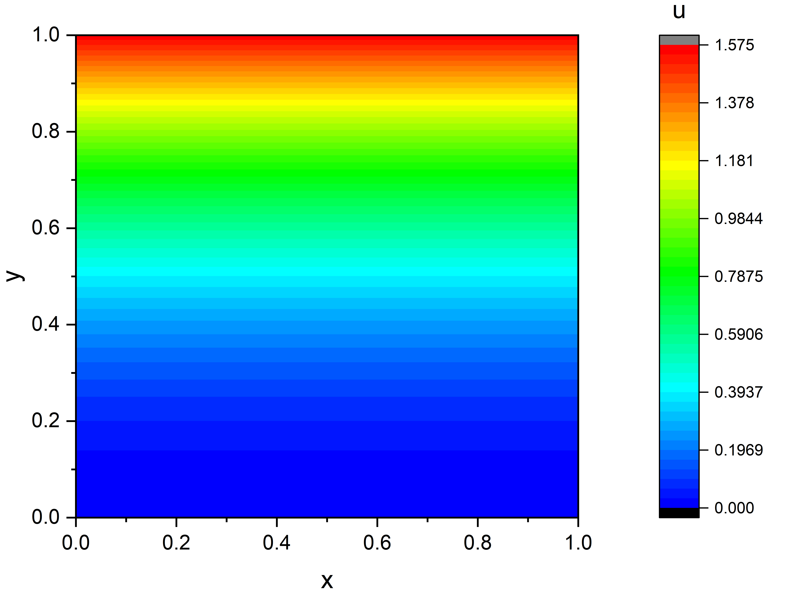}
 \end{subfigure}%
 \caption{Plot of finite element solution with 32 elements}
 \label{fig_con_square}
\end{figure}

The plot of the numerical solution is given in Figure~\ref{fig_con_square}.

The following table indicates the convergence of the numerical solution for each of the computational meshes considered. 
\begin{table}[H]
\centering
\caption{Different error measures for cubic-order triangular elements}
\begin{tabular}{|p{0.7 in}|p{1.25 in}|p{0.8 in}|p{1.25 in}|p{1.25 in}|}\hline
Elements & $E_{a}$ & $E_{r}$ & $l_{2}$-error norm\\ \hline
8 &  $7.148\times 10 ^{-8}$ &  0.000019 &  $1.2230 \times 10^{-7}$ \\ \hline
16 &  $5.190\times 10 ^{-8}$ & 0.000023 &  $1.1877 \times 10^{-7}$ \\ \hline
32 &  $4.800\times 10 ^{-8}$ & 0.000019 &  $1.1442 \times 10^{-7}$\\ \hline
\end{tabular}
 \label{ex1_res}
\end{table}
Here the symbols $E_{a}$ and $E_{r}$ are the absolute and relative errors, respectively. 
It is clear from Table~\ref{ex1_res} that our proposed finite element algorithm coupled with Picard's type linearization is an excellent choice to approximate the solution to the class of quasilinear partial differential equation.

\subsection{Domain with V-notch}
\begin{figure}[H]
 \centering
 \noindent\begin{subfigure}[b]{0.4\textwidth}
  \includegraphics[width=1\textwidth]{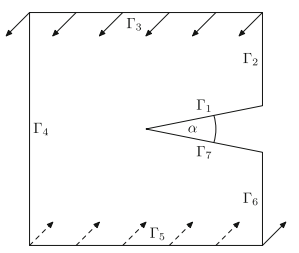}
 \end{subfigure}%
 \caption{V-notch Domain}
 \label{fig_vnotch}
\end{figure}

In this example, our goal is to present an efficient numerical method for the discretization of the quasilinear boundary value problem on a domain containing V-notch. For this setup, we do not have an exact analytical solution. The geometry and the boundary indicators are shown in Figure~\ref{fig_vnotch}. The boundary conditions are given in the Table~\ref{ex2_bc}.

\begin{table}[H]
\centering
\caption{Boundary Conditions for V-notch domain}
\begin{tabular}{p{0.9 in} p{0.5 in}}\hline
 Boundary & Values \\ \hline
 $\Gamma_{1}$ & 0 \\
 $\Gamma_{2}$ & 0 \\
 $\Gamma_{3}$ & $1-x$ \\
 $\Gamma_{4}$ & 1 \\
 $\Gamma_{5}$ & $1-x$ \\
 $\Gamma_{6}$ & $0$ \\
 $\Gamma_{7}$ & $0$ \\
 
\end{tabular}
 \label{ex2_bc}
\end{table}

\begin{figure}[H]
\centering 
  \includegraphics[width=0.4\textwidth]{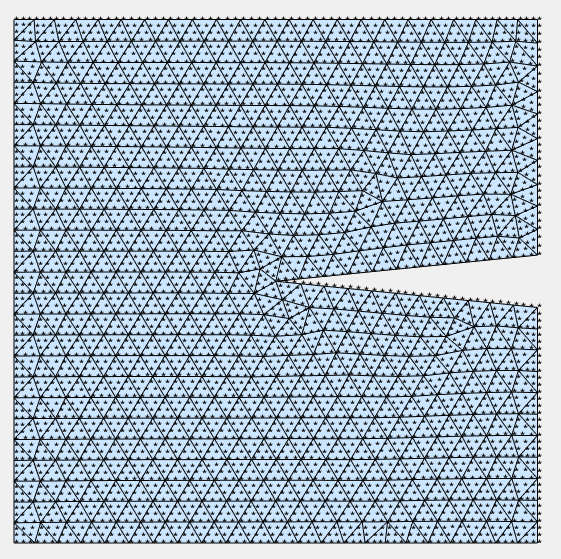}
 \caption{The computational mesh for the v-notch domain.}
 \label{fig_dis_vnotch}
\end{figure}

The computational mesh is depicted in Figure~\ref{fig_dis_vnotch} and the corresponding information about the elements, degrees-of-freedom, and the number of boundary nodes is given in Table~\ref{tab_vnotch}.  The mesh size in our computational discretization is taken as $0.045$. 

\begin{table}[H]
\centering
\caption{Details of domain discretization.  }
\begin{tabular}{p{0.6 in} p{0.5 in} p{1.25 in}}\hline
Elements & DOF & Boundary Nodes\\ \hline
1038 & 4843 & 342 \\ \hline
\end{tabular}
\label{tab_vnotch}
\end{table}

To monitor the convergence of our proposed method, we print out the maximum difference between the solutions at consecutive iterations.  The computed values are given in the Table~\ref{tab_vnotch}.  For a consistent comparison, we have computed the relative difference at a particular node in the mesh.  

\begin{table}[H]
\centering
\caption{The maximum value of the absolute difference between the solutions at two consecutive iterations. }
\begin{tabular}{|c|}\hline
 $ \max  | \Phi^{n}-\Phi^{n-1} |$ \\ \hline
0.079887041221999 \\
0.041362782097059 \\
0.021305707435756 \\
0.011408306865723 \\
0.006567123752435 \\
0.003989228093075 \\
0.002478142679311 \\
0.001569657000410 \\
0.001011378448757 \\
0.000661660540590 \\
0.000438859611370 \\
0.000294760677603 \\
0.000200279446002 \\
0.000137545640890 \\
0.000095400576785 \\
0.000066775278713 \\
0.000047132288683 \\
0.000033523256409 \\
0.000024010189043 \\
0.000017305167042 \\
0.000012543266832 \\
0.000009137782930 \\ \hline
\end{tabular}
 \label{tab_vnotch}
\end{table}

The computed solution at the end of the Picard's iteration is shown below in Figure~\ref{fig_con_vnotch}. 

\begin{figure}[H]
 \centering
  \includegraphics[width=0.5\textwidth]{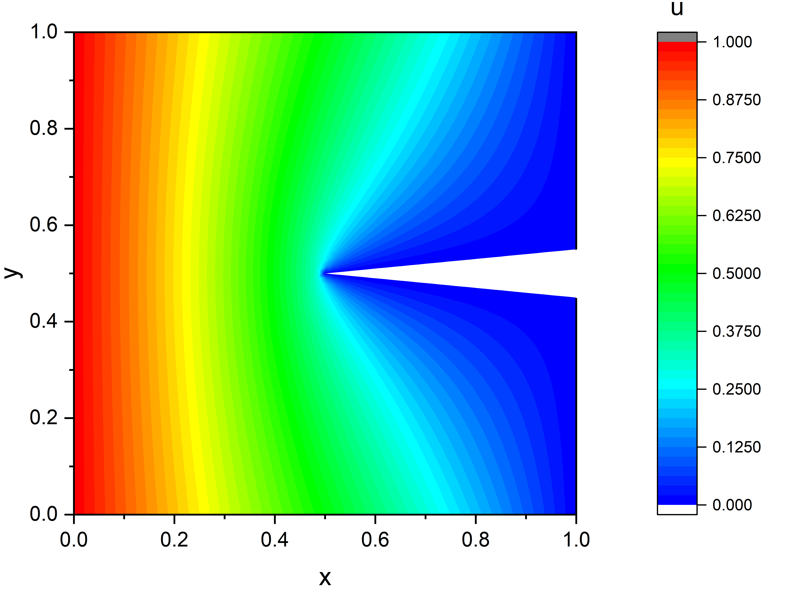}
 \caption{Computed solution from the finite element algorithm. }
 \label{fig_con_vnotch}
\end{figure}

It is clear from Table~\ref{tab_vnotch} that our proposed method takes less than $25$ iterations to obtain a reasonably accurate numerical solution. 

\subsection{Domain with V-notch and inclusions}
\begin{figure}[H]
 \centering 
  \includegraphics[width=0.4\textwidth]{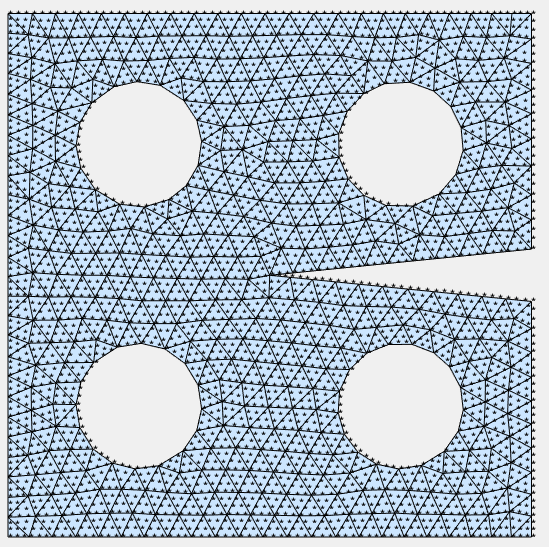}
 \caption{Traingulation of a material body containing both v-notch and inclusions. }
 \label{fig_vnotchwtinclusion}
\end{figure}

Here, our main objective is to present an efficient finite element method to handle the nonlinear partial differential equation posed on a domain containing v-notch and circular inclusions (or holes). The boundary conditions are kept the same as in the last example, and the inclusions are kept traction-free. The challenge for any numerical method is to handle the nonlinearities along with the triangulation of the domain. The proposed method studied in this paper considers the region near the inclusions as circular contrary to the other popular triangulations that use the points on the boundary of the inclusions to create the mesh with straight edges. The apparent bottleneck in many other finite element implementations is creating stable (good) computational mesh to efficiently evaluate the integrals needed to assemble the matrices and vectors. We do not have these issues as we do curved triangular elements along the circular inclusions. The triangulation of the domain is shown in Figure~\ref{fig_vnotchwtinclusion}.

The statistics of the mesh are given in the following Table~\ref{tab_vnotchwtincl}.

\begin{table}[H]
\centering
\caption{Details of the computational mesh. }
\begin{tabular}{p{0.6 in} p{0.5 in} p{1.25 in}}\hline
Elements & DOF & Boundary Nodes\\ \hline
836 & 4026 & 339 \\ \hline
\end{tabular}
\label{tab_vnotchwtincl}
\end{table}

To monitor the convergence of the numerical solution, we compute the maximum difference between two consecutive solutions and the same is depicted in Table~\ref{tab_vnotchwtincl}.

\begin{table}[H]
\centering
\caption{The maximum value of absolute difference between two consecutive iterations are tabulated below for V-notch domain}
\begin{tabular}{|c|}\hline
  $\max  | \Phi^{n}-\Phi^{n-1}|$ \\ \hline
0.061362195377219	\\
0.034003719699686	\\
0.019375215579367	\\
0.011322928480591	\\
0.006773972907949	\\
0.004099274297780	\\
0.002509970577362	\\
0.001554095858205	\\
0.000983314280567	\\
0.000642176317155	\\
0.000423204575014	\\
0.000281269726815	\\
0.000188420856898	\\
0.000127158649975	\\
0.000086411758823	\\
0.000059105537844	\\
0.000040676772821	\\
0.000028156283765	\\
0.000019596217704	\\
0.000013708893668	\\
0.000009636885800	\\ \hline
\end{tabular}
\label{tab_vnotchwtincl}
\end{table}

\begin{figure}[H]
 \centering
  \includegraphics[width=0.4\textwidth]{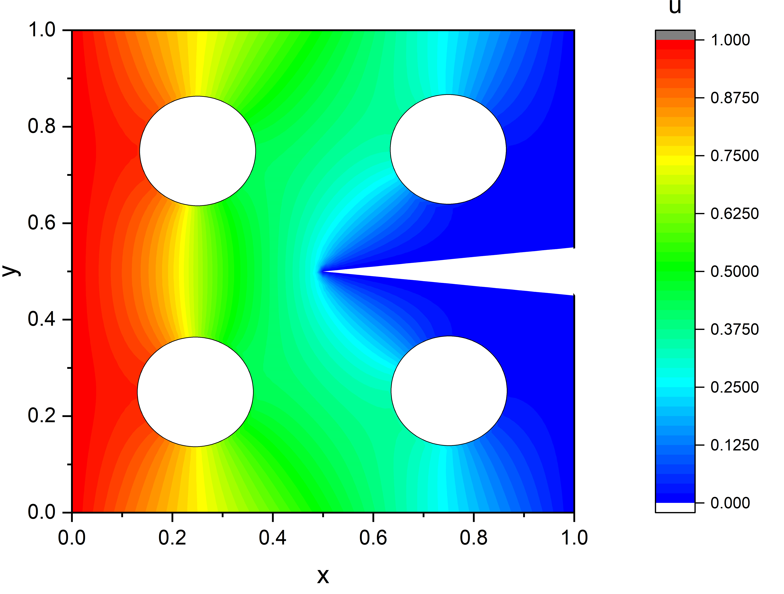}
 \caption{Plot of the numerical solution. }
 \label{fig_con_vnotchwtinclusions}
\end{figure}

Figure~\ref{fig_con_vnotchwtinclusions} shows the plot of the finite element solution for a strain-limiting material body containing both v-notch and circular inclusions. The solution presented is the one obtained after $21$ Picard's iteration and the overall algorithm met the stopping criteria. 

\section{Conclusion}\label{conclusion}
An efficient finite element for simulating the response of the algebraically nonlinear and geometrically linear strain-limiting elastic material containing both v-notch and circular inclusions is presented. A boundary value problem consists of a second-order quasilinear partial differential equation obtained by considering a special constitutive relationship within the setting of a new theory formulated in \cite{rajagopal2007elasticity,rajagopal2011modeling}. There are no analytical methods available to obtain a closed-form solution to such a setup. Hence, a novel type of finite element is proposed in this work. A general curved triangular element which is of cubic-order isoparametric type is utilized. Such curved elements are used along the circumference of the inclusions and which is adaptable to the analogous with the other triangular mesh elements in the bulk material body. A main element of our method is in approximating the curved boundary by second-order polynomials for the cubic order triangular method and also in constructing the efficient numerical method to approximate the primitive displacement variable. The proposed method is applied to three different boundary value problems. The accuracy of the numerical solution is tested against a manufactured solution.  The numerical results obtained are very good and improve the accuracy of the approximation of the quasi-linear partial differential equation. The standard Picard's type linearization takes less iteration to converge to the solution. Overall, the method presented in this paper is very good and crucial for several applications.  The approach used to simulate the elastic material under mechanical loading can be easily adaptable to the quasi-static crack propagation both in elastic materials whose response is described by algebraically nonlinear relationship \cite{yoon2021quasi,lee2022finite} and in materials whose parameters dependent upon the density \cite{HCYSMMDDB2024}, multiscale situations \cite{MVSMM2023,vasilyeva2023}, and multi-physics problems \cite{SMMDDB2023}.

\bibliographystyle{plain}
\bibliography{vnotchinclusion_reference}
 
\end{document}